\documentclass[12pt ]{article}
\usepackage{mathrsfs}
\usepackage{amsfonts}
\usepackage{amsfonts}
\usepackage{amsmath}
\usepackage{amssymb}
\usepackage{amsthm}
\newcommand{\bc}{\begin{center}}
\newcommand{\ec}{\end{center}}
\newcommand{\ba}{\begin{array}}
\newcommand{\ea}{\end{array}}

\date{}

\newtheorem{theorem}{Theorem}[section]
\newtheorem{corollary}{Corollary}[section]

\title{ An Elliptic Type Gradient Estimate For  the  Schr\"{o}dinger Equation }
\author{ Qihua Ruan}
\begin{document}

\maketitle
\begin {center}
{ \footnotesize School of Mathematics , Zhongshan University,
Guangzhou, Guangdong, 510275, P.R.China;

  \footnotesize Department of Mathematics, Putian University, Putian, Fujian,
  351100, P.R.China.

\footnotesize Email:\,\,ruanqihua@163.com

}
\end{center}
\thispagestyle{empty}

\begin{abstract}  
In this paper, the author discusses the elliptic type gradient
estimate for the solution of the time-dependent Schr\"{o}dinger
equations on noncompact manifolds. As its application,  the
dimension-free Harnack inequality and   the Liouville type theorem
for the Schr\"{o}dinger equation  are proved.
 \vspace{0.3cm}


\vspace{0.3cm} \noindent {\bf Mathematics Subject Classification
(2000):}\ 58G11, 53C21.
\end{abstract}

\setcounter{equation}{0}
\section{Introduction}
In \cite{CY}, Cheng and Yau proved the classical gradient estimate
for harmonic functions.  Later Li and Yau \cite{LY} obtained the
following parabolic type gradient estimate for the heat equation
$u_{t}=\triangle u$ and where the subscript $t$ denotes the partial
differentiation with respect to the $t$-variable, $\triangle$ is the
Laplacian operator on $M$.
\begin{theorem}$(Li-Yau)$ Let $M$ be an n-dimensional complete manifold with
 Ricci curvature  bounded below by $-K$,
$K\geq0$. Suppose that $u$ is any positive solution to the heat
equation $u_{t}=\triangle u$ in  $B(x_{0}, 2R)\times[t_{0}-2T,
t_{0}]$, then for $\forall\alpha>1$,
\begin{eqnarray}\label{1.1}
\frac{|\nabla u|^{2}}{u^{2}}-\alpha\frac{u_{t}}{u}\leq
\frac{C}{R^{2}}+\frac{n\alpha^{2}}{2T}
+\frac{n\alpha^{2}}{\sqrt{2}(\alpha-1)}K
\end{eqnarray}
in $B(x_{0}, R)\times[t_{0}-T, t_{0}]$. Here $\nabla$ denotes the
gradient operator on $M$, and the positive constant $C$ depends only
on dimension $n$.
\end{theorem}
Recently X.D.Li \cite{Li} studied the heat equation $u_{t}=\triangle
u+\nabla \phi\cdot\nabla u$, $\phi\in C^{2}(M)$,  on  a manifold
with Bakry-Emery's Ricci curvature
$\widetilde{Ric}:=Ric-\nabla^{2}\phi-\frac{1}{m-n}\nabla\phi\otimes\nabla\phi$,
where the constant $m> n$,  $Ric$ denotes the Ricci curvature of $M$
and $m=n$ if and only if $\phi=0$. If the curvature condition in
Theorem 1.1 is replaced by the following curvature condition:
\begin{eqnarray}\label{1.2}
\widetilde{Ric}\geq-K,
\end{eqnarray}
then Li obtained a similar gradient estimate for the heat equation
$u_{t}=\triangle u+\nabla \phi\cdot\nabla u$: $$\frac{|\nabla
u|^{2}}{u^{2}}-\alpha\frac{u_{t}}{u}\leq
\frac{C}{R^{2}}+\frac{m\alpha^{2}}{2T}
+\frac{m\alpha^{2}}{\sqrt{2}(\alpha-1)}K.$$

 The Bakry-Emery Ricci
curvature has some interesting properties. In \cite{R1}, the author
established some connections between the Bakry-Emery Ricci curvature
and Ricci flow. In \cite{R2}, he also proved some Liouville type
theorems for the above heat equation $u_{t}=\triangle u+\nabla
\phi\cdot\nabla u$ on a Manifold with the Bakry-Emery Ricci
curvature.

As was shown in \cite{LY},  using (1.1) we can derive the following
Harnack inequality:
$$u(x_{1}, t_{1})\leq u(x_{2},
t_{2})(\frac{t_{2}}{t_{1}})^{\frac{n\alpha}{2}}
\exp(\frac{\alpha\rho^{2}(x_{1},
x_{2})}{4(t_{2}-t_{1})}+\frac{n\alpha
K}{\sqrt{2}(\alpha-1)}(t_{2}-t_{1}))$$ where  $\forall x_{1},
x_{2}\in M$, $\rho(x_{1}, x_{2})$ denotes the geodesic distance
between $x_{1}$ and $x_{2}$, and  $0<t_{1}<t_{2}<+\infty$. We notice
that from this kind of Harnack inequality we can only compare the
solutions  at different times.  In \cite{H}, Hamilton got the
following elliptic type gradient estimate on compact manifolds. From
this gradient estimate one can compare the solution of two different
points at the same time.
\begin{theorem}$(Hamilton)$ Let M be a compact manifold without
boundary and with Ricci curvature  bounded below by $-K$, $K\geq0$.
Suppose that $u$ is any positive solution to the heat equation
$u_{t}=\triangle u$ with $u\leq C$ for all $(x, t)\in M\times(0,
+\infty)$.  Then
\begin{eqnarray}\label{1.3}
\frac{|\nabla u|^{2}}{u^{2}}\leq(\frac{1}{t}+2K)(\ln\frac{C}{u}).
\end{eqnarray}
\end{theorem}
  Recently Souplet and Zhang \cite{SZ} generalized the
elliptic type gradient estimate to noncompact manifolds.
\begin{theorem}$(Souplet-Zhang )$ Let $M$ be an n-dimensional complete noncompact manifold with
 Ricci curvature  bounded below by $-K$,
$K\geq0$. Suppose that $u$ is any positive solution to the heat
equation $u_{t}=\triangle u$ in  $Q_{2R, 2T}\equiv B(x_{0},
2R)\times[t_{0}-2T, t_{0}]$,  and $u\leq C$ in $Q_{2R, 2T}$.  Then
\begin{eqnarray}\label{1.4}
\frac{|\nabla u|}{u}\leq C_{1}(
\frac{1}{R}+\frac{1}{T^{\frac{1}{2}}} +\sqrt{K})(1+\ln\frac{C}{u})
\end{eqnarray} in $Q_{R, T}$. Here  the
positive constant $C_{1}$  dependes only on dimension $n$, moreover
$C_{1}$ increases to infinity as  $n$ goes to infinity.
\end{theorem}
We know that (1.4) is a local gradient estimate. Letting
$R\rightarrow+\infty$, we can obtain the global gradient estimate:
\begin{eqnarray}\label{1.5}
\frac{|\nabla u|}{u}\leq C_{1}( \frac{1}{T^{\frac{1}{2}}}
+\sqrt{K})(1+\ln\frac{C}{u})
\end{eqnarray} We notice that the
constant $C_{1}$ in (1.5) depends on the dimension $n$, moreover
 $C_{1}$ increases to infinity as  $n$ goes to infinity,  so this
gradient estimate cannot be applied to infinite dimensional
manifolds. However Hamilton's gradient estimate (1.3) does not
depend on dimension $n$. In this paper, we want to prove the
following dimension-free elliptic type gradient estimate for a more
general equation---the Schr\"{o}dinger equation with potential $h(x,
t)$:
\begin{eqnarray}\label{1.6}u_{t}=\triangle u+\nabla\phi\cdot\nabla
u+hu.
\end{eqnarray}

\begin{theorem}Let $M$ be an n-dimensional  complete noncompact
manifold satisfying the curvature condition (1.2). Suppose that the
potential $h$ is a negative function defined on  $M\times(0,
+\infty)$ which is $C^{1}$ in the $x$-variable,  and that $u$ is any
positive solution to the Schr\"{o}dinger equation (1.6) with $u\leq
C$ for all $(x, t)\in M\times(0, +\infty)$. Then
\begin{eqnarray}\label{1.7}
\frac{|\nabla u|}{u}\leq ( \frac{1}{t^{\frac{1}{2}}}
+\sqrt{2K}+|\nabla
\sqrt{-h}|^{\frac{1}{2}})(1+\ln\frac{C}{u}).\end{eqnarray}
\end{theorem}

As proved by Li and Yau in \cite{LY}, we can derive the
dimension-free Harnack inequality for the equation (1.8). We must
point out that the dimension-free Harnack inequality for the heat
equation $u_{t}=\triangle u+\nabla\phi\cdot\nabla u$ was first
observed by Wang \cite{W1}. This kind of inequality can be widely
applied in geometric analysis and probability, more details can be
seen in Wang's survey article\cite{W2}.
\begin{corollary}Let $M$ be an n-dimensional  complete noncompact
manifold satisfying the curvature condition (1.2). Suppose that the
potential $h$ is a negative function defined on  $M\times(0,
+\infty)$ which is $C^{1}$ in the $x$-variable, and $|\nabla
\sqrt{-h}|\leq C_{2}$  for a positive constant $C_{2}$,  and that
$u$ is any positive solution to the Schr\"{o}dinger equation (1.6)
with $u\leq 1$ for all $(x, t)\in M\times(0, +\infty)$.  Then for
$\forall x_{1}, x_{2}\in M$,

$$u(x_{2}, t)\leq u(x_{1}, t)^{\beta}e^{1-\beta}$$
where
$\beta=\exp(-\frac{\rho}{t^{\frac{1}{2}}}-(\sqrt{2K}+\sqrt{C_{2}})\rho)$
and
 $\rho=\rho(x_{1},
x_{2})$ denotes the geodesic distance between $x_{1}$ and $x_{2}$.
\end{corollary}

When potential $h(x, t)=-\lambda$, $\lambda$ is a positive constant,
(1.6) becomes
\begin{eqnarray}\label{1.8}
u_{t}=\triangle u+\nabla\phi\cdot\nabla u-\lambda u.
\end{eqnarray}Then we can easily obtain the following
dimension-free elliptic type gradient estimate for the equation
(1.8), which improves Souplet and Zhang's gradient estimate (1.5).
\begin{corollary}Let $M$ be an n-dimensional  complete noncompact
manifold satisfying the curvature condition (1.2). Suppose that $u$
is any positive solution to the equation (1.8) with $u\leq C$ for
any $(x, t)\in M\times(0, +\infty)$.  Then
\begin{eqnarray}\label{1.9}
\frac{|\nabla u|}{u}\leq ( \frac{1}{t^{\frac{1}{2}}}
+\sqrt{2K})(1+\ln\frac{C}{u}).
\end{eqnarray}
\end{corollary}

We notice that when the solution of (1.8) does not depend on time,
the equation (1.8) becomes
\begin{eqnarray}\label{1.10}
\triangle u(x)+\nabla\phi(x)\cdot\nabla u(x)-\lambda u(x)=0.
\end{eqnarray} From
(1.9), letting $t\rightarrow\infty$, we have the gradient estimate
for the equation (1.10):
\begin{eqnarray}\label{1.11}
\frac{|\nabla u|}{u}\leq \sqrt{2K}(1+\ln\frac{C}{u}).
\end{eqnarray}
 From Schoen and Yau's nice
book \cite{SY}, it is known that the gradient estimate of the
positive solution for the equation $\triangle u-\lambda u=0$ is that
$\frac{|\nabla u|}{u}\leq C(n, K, \lambda)$, where the positive
constant $C(n, K, \lambda)$ depends on $n,\ K,\ \lambda$. However
the gradient estimate (1.11) doesn't depend on $\lambda$, so this
result allows us to deduce the following Liouville type theorem:
\begin{corollary}Let $M$ be an n-dimensional  complete noncompact
manifold with  nonnegative Bakry-Emery's Ricci curvature, i.e. $K=0$
in (1.2). Then there does not exist any positive and bounded
solution to Schr\"{o}dinger equation (1.10).
\end{corollary}

{\bf Remark:} The Liouville type theorem for the Schr\"{o}dinger
equation was investigated by many authors. In \cite{L}, \cite{LY},
\cite{N}
 and \cite{M}, the authors considered the Schr\"{o}dinger equation
with nonnegative potential and obtained a similar Liouville type
theorem. However in our case the potential is negative. In
\cite{G1}, \cite{GH} and \cite{KL}, the authors discussed the
Liouville type theorem for the Schr\"{o}dinger equation on parabolic
manifolds. However in our case,  the manifold need not be parabolic.

\setcounter{equation}{0}
\section{Proof of Main  Theorem and Corollaries  }
In this section, we will prove Theorem 1.4 and  its corollaries.

{\bf Proof:}First we define the operators $L=\triangle+\nabla\phi$
and $\Box =L-\partial_{t}$. Since $u$ is the solution of the
equation $\Box u=-hu$, we notice that the solution $u$ is invariant
under scaling, so we can assume $u\leq1$. Setting $f=\ln u$,
$w=|\nabla \ln (1-f)|^{2}$, then $f\leq0$ and also satisfies the
following equation
\begin{eqnarray}\label{2.1}
\Box f= u^{-1}\Box u-|\nabla f|^{2}=-h-|\nabla f|^{2}.
\end{eqnarray}
By (2.1), we see that
\begin{eqnarray}\label{2.2}
\Box\ln(1-f)=-\frac{\Box f}{1-f}-|\nabla \ln
(1-f)|^{2}=\frac{h}{1-f}-f w.
\end{eqnarray}

For convenience of computation, we introduce the two curvature
operators defined by  Bakry-Emery \cite{BE}.  We define the
curvature operator $\Gamma$ by
$$\Gamma (f, g)=\frac{1}{2}\{L(fg)-fLg-gLf\}, \
\forall f, g \in C^{2}(M)$$

It is obvious that \begin{eqnarray}\label{2.3}\Gamma(f,f)=|\nabla
f|^{2} \end{eqnarray}

And the Bakry-Emery curvature operator $\Gamma_{2}$ is defined to be
a bilinear map: $$\Gamma_{2}(f, g)=\frac{1}{2}\{L\Gamma(f,
g)-\Gamma(Lf, g)-\Gamma(f, Lg)\}, \  \forall f, g \in C^{2}(M)$$ In
particular, we have that $$
\Gamma_{2}(f,f)=\frac{1}{2}\{L\Gamma(f,f)-2\Gamma(L f,f)\}$$ Through
direct calculation, we find that the Bakry-Emery curvature operator
$\Gamma_{2}$ also satisfies
\begin{eqnarray}\label{2.4}
\Gamma_{2}(f,f)=\frac{1}{2}\{\Box\Gamma(f,f)-2\Gamma(\Box f,f)\},\
\forall f\in C^{2}(M).\end{eqnarray} In the sequel, for convenience
we sometimes denote $\Gamma(f,f)$ and $\Gamma_{2}(f,f)$ by
$\Gamma(f)$ and $\Gamma_{2}(f)$ respectively.

By Bochner's formula, under the curvature condition (1.2), we have
that
 \begin{eqnarray}\label{2.5}\Gamma_{2}(f) \geq-K\Gamma f,\
\forall f\in C^{2}(M).\end{eqnarray} This inequality was proved  in
\cite{BE} (see also  \cite{Li}).

 By $(2.2)-(2.5)$, we know that
 $$\Box w=\Box\Gamma(\ln(1-f))\hspace{5.1cm}$$
 $$=2\Gamma_{2}(\ln(1-f))+2\Gamma(\Box\ln(1-f), \ln(1-f))$$
 $$\hspace{2.4cm}=2\Gamma_{2}(\ln(1-f))+2\Gamma(\frac{h}{1-f}, \ln(1-f))-2\Gamma(fw, \ln(1-f))$$
 $$\hspace{2.4cm}\geq -2Kw-\frac{2h}{1-f}w-\frac{2|\nabla h|}{1-f}w^{\frac{1}{2}}+2(1-f)w^{2}+2f|\nabla w|w^{\frac{1}{2}}$$
$$\hspace{0.8cm}\geq -2Kw-\frac{|\nabla\sqrt{-h}|^{2}}{1-f}+2(1-f)w^{2}+2f|\nabla w|w^{\frac{1}{2}}$$
 $$\hspace{0.8cm}\geq -2Kw-|\nabla \sqrt{-h}|^{2}+2(1-f)w^{2}+2f|\nabla w|w^{\frac{1}{2}}$$
where we use the fact that $h<0$, $\frac{1}{1-f}<1$ in the above
last two inequalities. Hence  we obtain
\begin{eqnarray}\label{2.6}
\Box w\geq -2Kw-|\nabla \sqrt{-h}|^{2}+2(1-f)w^{2}+2f|\nabla
w|w^{\frac{1}{2}}.\end{eqnarray}

Now we choose a cut-off function $\eta$ to be a $C^{2}$ function on
$[0, +\infty)$ satisfying \[ \eta(t)=1, \qquad\mbox{for}\qquad\
0\leq t\leq1; \qquad \eta(t)=0,\qquad\mbox {for}\qquad  t\geq2;
\] and
\[0\leq\eta(t)\leq1,\ \ -C\eta(t)^{\frac{1}{2}}\leq\eta'(t)\leq0,  \
\ \eta''(t)\geq-C\qquad\mbox{for}\qquad\   \forall t\geq0.\] Let
$\rho(p, x)$ be the geodesic distance between $p$ and $x$, and
define $\psi(x)=\eta(\frac{\rho(p, x)}{R})$. Then we have that
\begin{eqnarray}\label{2.7}\frac{|\nabla\psi|^{2}}{\psi}=
\frac{|\eta'|^{2}|\nabla
\rho|^{2}}{R^{2}\eta}=\frac{|\eta'|^{2}}{R^{2}\eta}\leq\frac{C^{2}}{R^{2}}.\end{eqnarray}
and
\begin{eqnarray}\label{2.8}L\psi(x)=\frac{\eta''|\nabla \rho|^{2}}{R^{2}}+\frac{\eta'L \rho}{R}\geq-\frac{mC}{R^{2}}
-\frac{mK}{R},\  x\notin cut(p),\end{eqnarray} where we use a
general Laplacian comparison theorem: $L
\rho\leq\frac{m-1}{\rho}+(m-1)K$ under the curvature condition
(1.2), which was proved by Qian \cite{Q}.

Let $\varphi=t\psi$. Suppose that $\varphi w$ attains its maximum at
the point $(x_{0}, t_{0} )\in B(p, 2R)\times[0, T]$. According a
well known argument of Calabi\cite{C},  we can assume that $x_{0}$
is not in the cut locus of $p$. Then at $(x_{0}, t_{0} )$ we have
$$\nabla(\varphi w)=0,\ \triangle(\varphi w)\leq0,\ \frac{\partial}{\partial t}(\varphi w)\geq0$$
Thus we derive that at $(x_{0}, t_{0} )$
\begin{eqnarray}\label{2.9}\Box(\varphi w)\leq0,\ \nabla w=-\frac{\nabla
\varphi}{\varphi}w.
\end{eqnarray}
Hence (2.9) implies that at $(x_{0}, t_{0} )$
$$\varphi\Box w+\Box \varphi\omega+2\Gamma(\varphi,
w)\leq0$$Combining the above inequality with (2.6) and (2.9),
 we
have at $(x_{0}, t_{0} )$
$$\varphi\{-2Kw-|\nabla \sqrt{-h}|^{2}+2(1-f)w^{2}+2f\frac{|\nabla\varphi|}{\varphi}w^{\frac{3}{2}}\}+\Box \varphi\omega
-\frac{2|\nabla\varphi|^{2}}{\varphi}w\leq0$$ By the Young
inequality: $2ab\leq a^{2}+b^{2}$ and $f<0$,  we see at $(x_{0},
t_{0} )$
$$-2K\varphi w-|\nabla \sqrt{-h}|^{2}\varphi+(1-f)\varphi w^{2}-\frac{f^{2}|\nabla\varphi|^{2}}{(1-f)\varphi}w
+\Box \varphi\omega-\frac{2|\nabla\varphi|^{2}}{\varphi}w\leq0$$
Multiplying the  both sides of the above inequality by
$\frac{\varphi}{1-f}$, we notice  that $0\leq\psi\leq1$,
$\frac{f^{2}}{(1-f)^{2}}\leq1$ and $\frac{1}{1-f}\leq1$, then from
(2.7) and (2.8) we know at $(x_{0}, t_{0} )$
$$(\varphi w)^{2}-T(2K+\frac{3C^{2}}{R^{2}}+\frac{mC}{R^{2}}
+\frac{mK}{R}+\frac{1}{T})(\varphi w)-|\nabla
\sqrt{-h}|^{2}T^{2}\leq0$$ Since $\psi=1$ on $B(p, R)$,  then we
have
$$w(x, T)\leq 2K+\frac{3C^{2}}{R^{2}}+\frac{mC}{R^{2}}
+\frac{mK}{R}+\frac{1}{T}+|\nabla \sqrt{-h}|$$ for $\forall x\in
B(p, R)$. Since $T$ is arbitrary, let $R\rightarrow +\infty$, then
we obtain
$$\frac{|\nabla u|}{u}\leq ( \frac{1}{t^{\frac{1}{2}}}
+\sqrt{2K}+|\nabla \sqrt{-h}|^{\frac{1}{2}})(1-\ln u)$$ The scaling
invariance of $u$ implies (1.7). This completes the proof of Theorem
1.4.

Now we give the proof of Corollary 1.5:

{\bf Proof:} Let $\gamma(s)$ be a minimal geodesic joining $x_{1}$
and $x_{2}$, $\gamma:[0, 1]\rightarrow M$, $\gamma(0)=x_{2}$,
$\gamma(1)=x_{1}$.  Then by (1.9)
$$\ln\frac{1-f(x_{1}, t)}{1-f(x_{2}, t)}=\int_{0}^{1}\frac{d\ln(1-f(\gamma(s), t))}{ds}ds$$
$$\hspace{1.9cm}\leq\int_{0}^{1}|\dot{\gamma}|\cdot\frac{|\nabla u|}{u(1-\ln u)}ds$$
$$\hspace{2.1cm}\leq\frac{\rho}{t^{\frac{1}{2}}}+(\sqrt{2K}+\sqrt{C_{2}})\rho$$
Then
\begin{eqnarray}\label{2.10}\frac{1-f(x_{1}, t)}{1-f(x_{2}, t)}\leq \exp\{\frac{\rho}{t^{\frac{1}{2}}}
+(\sqrt{2K}+\sqrt{C_{2}})\rho\}.
\end{eqnarray} Let
$\beta=\exp(-\frac{\rho}{t^{\frac{1}{2}}}-(\sqrt{2K}+\sqrt{C_{2}})\rho)$,
(2.10) then implies that
$$\frac{1-\ln u(x_{1}, t)}{1-\ln u(x_{2}, t)}\leq \frac{1}{\beta}$$
Through some easy computation, we obtain
$$u(x_{2}, t)\leq u(x_{1}, t)^{\beta}e^{1-\beta}$$
This completes the proof of Corollary 1.5.

Since the potential in (1.8) is a constant,   the term
$|\nabla\sqrt{-h}|^{\frac{1}{2}}$ in (1.7) is equal to zero.  So
 Corollary 1.6 follows by (1.7).   Let $t\rightarrow \infty$ and $K=0$ in (1.9),
then $|\nabla \log u|=0$. Thus $u$ must be a constant. From (1.10),
we know $u\equiv 0$. So there does not exist any positive and
bounded solution of the equation (1.10). This completes the
Corollary 1.7.

{\bf Acknowledgements:} The author wishes to express his thank to
Professor Z.H.Chen and Professor X.P.Zhu for their constant
encouragement.

\end{document}